\documentclass[11pt]{article}
\usepackage{mathrsfs}
\usepackage[centertags]{amsmath}
\usepackage{amsfonts}
\usepackage{amssymb}
\usepackage{amsthm}
\usepackage{cases}
\usepackage{indentfirst}
\usepackage{epsfig}
\usepackage {graphicx}
\usepackage{color}

\date{}

\definecolor{c20}{rgb}{0.,0.7,0.}
\definecolor{c30}{rgb}{0.,0.,1.}
\definecolor{c40}{rgb}{1,0.1,0.7}
\definecolor{c50}{rgb}{1,0,0}

\def\peng#1{\textcolor{c30}{#1}}
\def\peng#1{#1}

\newtheorem{theorem}{Theorem}
\newtheorem{corollary}{Corollary}
\newtheorem{remark}{Remark}
\numberwithin{equation}{section}

\def\P{\operatorname*{\mathbf{P}}}
\def\E{\operatorname*{\mathbf{E}}}

\def\R{\operatorname*{\mathbb{R}}}

\allowdisplaybreaks[4]
\usepackage{setspace}
\begin{document}
\title{\textbf{Maxima and minima of independent and  non-identically distributed bivariate Gaussian triangular arrays  }}
\author{{Yingyin Lu \;\;\; Zuoxiang Peng\thanks{Corresponding author. Email: pzx@swu.edu.cn} }    \\
{\small School of Mathematics and Statistics, Southwest University, Chongqing, 400715, China}}

\maketitle
\begin{quote}
{\bf Abstract.}~~In this paper, joint limit distributions of maxima and minima on independent and
non-identically distributed bivariate Gaussian triangular arrays is derived as the correlation coefficient of $i$th vector of given $n$th row is the function of $i/n$.
Furthermore, second-order expansions of joint distributions of maxima and minima are established \peng{if the correlation} function satisfies some regular conditions.

{\bf Keywords.}~~Bivariate Gaussian triangular arrays; maximium and minimium; limiting distribution; second-order expansion.
\end{quote}

\section{Introduction}
\label{sec1}

Let $\{\left(\xi_{ni},\eta_{ni}\right),1\leq i \leq n,n\geq 1\}$ be independent bivariate Gaussian triangular arrays with
$\E \xi_{ni}=\E \eta_{ni}=0$, $\E \xi_{ni}^{2}=\E \eta_{ni}^{2}=1$, and let $\rho_{ni}$ denote the correlation coefficient
of $\left(\xi_{ni},\eta_{ni}\right)$,$ 1\leq i \leq n $. The bivariate maxima $\mathbf{M}_{n}$ is defined  componentwise by
\[\mathbf{M}_{n}=(M_{n1},M_{n2})=\left(\max_{1 \leq i\leq n} \xi_{ni}, \max_{1\leq i \leq n} \eta_{ni}\right).\]
For the case of $\rho_{ni}=\rho_{n}$, the seminal paper of H\"{u}sler and Reiss (1989) \peng{showed that the limiting} distribution of normalized maxima of such bivariate Gaussian triangular arrays is
\begin{equation}\label{HRD}
\lim_{n\to\infty}\P\left(M_{n1}\le x/b_{n}+b_{n},M_{n2}\le y/b_{n}+b_{n}\right)=H_{\lambda}(x,y)
\end{equation}
provided that the following so-called H\"{u}sler-Reiss condition
\begin{eqnarray}\label{HRC}
\lim_{n\to\infty} b_{n}^2\left(1-\rho_{n}\right)=2\lambda^2 \quad\mbox{with}\; \lambda\in\left[0,\infty\right]
\end{eqnarray}
holds, where the norming constant $b_{n}$ satisfies
\begin{eqnarray}
\label{eq1.2}
1-\Phi(b_{n})=\frac{1}{n}
\end{eqnarray}
and
\begin{eqnarray*}
H_{\lambda}(x,y)= \exp\left(
-\Phi\left(\lambda+\frac{x-y}{2\lambda}\right)e^{-y} -
\Phi\left(\lambda+\frac{y-x}{2\lambda}\right)e^{-x}
\right)
\end{eqnarray*}
with $H_{0}(x,y)=\Lambda\left(\min(x,y)\right)$ and $H_{\infty}(x,y)=\Lambda(x)\Lambda(y)$ for all $ (x,y)\in{\R}^{2}$, \peng{where $\Lambda(x)=\exp(-e^{-x})$}, and $\Phi(x)$ denotes the standard Gaussian distribution. Kabluchko et al. (2009) showed that \eqref{HRC} also is the necessary condition for \eqref{HRD}.

Liao and Peng (2016 ) extended the work of H\"{u}sler and Reiss (1989) to independent and non-identically distributed bivariate Gaussian triangular arrays by assuming that the correlation $\rho_{ni}$ satisfying
\begin{eqnarray}
\label{eq1.4}
\rho_{ni}=1-\frac{m(i/n)}{\log n}, \quad 1\le i\le n
\end{eqnarray}
for some positive function $m(x)$ defined on $[0,1]$, and showed that
\begin{equation}\label{HRDC}
\lim_{n\to\infty}\P\left(M_{n1}\le x/b_{n}+b_{n},M_{n2}\le y/b_{n}+b_{n}\right)=H(x,y)
\end{equation}
if \eqref{eq1.4} holds, where
\begin{eqnarray}
\label{eq1.5}
H(x,y)=\exp\left(-e^{-y}\int_{0}^{1}\Phi\left(\sqrt{m(t)}+\frac{x-y}{2\sqrt{m(t)}}\right)dt
-e^{-x}\int_{0}^{1}\Phi\left(\sqrt{m(t)}+\frac{y-x}{2\sqrt{m(t)}}\right)dt\right).
\end{eqnarray}
\peng{Note that \eqref{HRDC} can be achieved by Theorem 2.5 in Engelke et al. (2015) by adding monotonicity to $m(x)$.}
Liao et al. (2016) considered the limiting distribution of  $(n(\max F_{1}(\xi_{ni})-1), n(\max F_{2}(\eta_{ni})-1))$ provided that each vector has Gaussian copula with $\rho_{ni}$ given by \eqref{eq1.4}, where $F_{1}$ and $F_{2}$ are distributions of $\xi_{ni}$ and $\eta_{ni}$, respectively. For more details, see Theorem 2.1 in Liao et al. (2016). For some other work related to \peng{H\"{u}sler-Reiss models} and its extensions, see, e.g., Hashorva (2005, 2006, 2013), Hashorva and Weng (2013), Hashorva et al. (2012), Hashorva et al. (2014), Frick and Reiss (2013) and D\c{e}bicki et al. (2014). Hsing et al. (1996) and French and Davis (2013) showed that Gaussian random fields with correlation between neighboring satisfying the conditions similar to \eqref{HRC} exhibit extremal clustering in the limits.

One interesting topic in extreme value theory is the convergence rates of distributions of order statistics to their ultimate extreme value distributions. For the univariate settings, this work was considered by
de Haan and Resnick (1996) under the second-order regular varying conditions, and Hall (1979) and Nair (1981) for independent and identically distributed Gaussian sequence. For bivariate settings under second-order regular varying conditions and other conditions, see de Haan and Peng (1997). Higher-order expansions and uniform convergence rates of joint distributions of maxima of bivariate Gaussian triangular arrays were derived respectively by Hashorva et al. (2016) and Liao and Peng (2014a) under the refined H\"{u}sler-Reiss conditions.

The objective of this paper is to establish the first and second-order asymptotics of the joint distributions of normalized
maxima and minima of the bivariate Gaussian triangular arrays with assumption \eqref{eq1.4}. Under the condition \eqref{HRC}, the joint asymptotics of maxima and minima of bivariate Gaussian
triangular arrays were studied by Liao and Peng (2014b). Precisely, let
\[\mathbf{m}_{n}=(m_{n1},m_{n2})=\left(\min_{1 \leq i\leq n} \xi_{ni}, \min_{1\leq i \leq n} \eta_{ni}\right)\]
denote the bivariate minima of the Gaussian triangular arrays, and
\begin{eqnarray}
\label{eq1.3}
\mathbf{v}_{n}=(-b_{n}+x_{1}/b_{n},
-b_{n}+y_{1}/b_{n})\quad\mbox{and}\quad
\mathbf{u}_{n}=(b_{n}+x_{2}/b_{n}, b_{n}+y_{2}/b_{n}),
\end{eqnarray}
Liao and Peng (2014b) showed that
\begin{eqnarray*}
\lim_{n\to\infty}\P\left(\mathbf{M}_{n}\leq \mathbf{u}_{n},\mathbf{m}_{n}\leq \mathbf{v}_{n}\right)=H_{\lambda}(x_{2},y_{2})\widetilde{H}_{\lambda}(x_{1},y_{1})
\end{eqnarray*}
if \eqref{HRC} holds, where
\[{\widetilde{H}_{\lambda}(x,y)=1-\Lambda(-x)-\Lambda(-y)+H_{\lambda}(-x,-y)}\]
with
\begin{equation}\label{minL}
\left\{{{\begin{array}{*{20}c}
{\widetilde{H}_{0}(x,y)=1-\Lambda(-x)-\Lambda(-y)+\Lambda\left(\min(-x,-y)\right)}\hfill\\
{\widetilde{H}_{\infty}(x,y)=1-\Lambda(-x)-\Lambda(-y)+\Lambda(-x)\Lambda(-y)}\hfill\\
\end{array}}}\right.
\end{equation}
for $(x,y)\in{\R}^{2}$. Furthermore, with $ \lambda_{n} = \left(\frac{1}{2} b_{n}^{2}(1-\rho_{n})\right)^\frac{1}{2} $, under the following
refined H\"{u}sler-Reiss condition
\begin{equation*}
\lim_{n\to\infty} b_{n}^{2}\left(\lambda_{n}-\lambda\right)=\alpha\; \mbox{with}\; \lambda\in (0,\infty),
\end{equation*}
Liao and Peng (2014b) showed that
\begin{eqnarray*}
& & \lim_{n\to \infty} b_{n}^{2}\left[ P(\mathbf{M}_{n}\leq \mathbf{u}_{n}, \mathbf{m}_{n} \leq \mathbf{v}_{n})
- H_{\lambda}(x_{2},y_{2})\tilde{H}_{\lambda}(x_{1},y_{1})  \right] \nonumber \\
&= &   H_{\lambda}(x_{2},y_{2})\Big[Q_{\lambda}(x_{2},y_{2}) +
H_{\lambda}(-x_{1},-y_{1}) \Big(Q_{\lambda}(-x_{1},-y_{1}) +
Q_{\lambda}(x_{2},y_{2})\Big)
\nonumber \\
&& \qquad\qquad\quad- \Lambda(-x_{1})\Big( \kappa(-x_{1}) +
Q_{\lambda}(x_{2},y_{2}) \Big) - \Lambda(-y_{1})\Big( \kappa(-y_{1})
+ Q_{\lambda}(x_{2},y_{2})  \Big)\Big].
\end{eqnarray*}
 Here,
\[Q_{\lambda}(x,y)= \kappa(x)\Phi\left( \lambda + \frac{y-x}{2\lambda} \right)
+ \kappa(y) \Phi\left( \lambda + \frac{x-y}{2\lambda} \right)
-\left( \lambda^{3} + \lambda x +\lambda y +2\lambda +2\alpha
\right)e^{-x}\varphi\left( \lambda + \frac{y-x}{2\lambda} \right),
\] where $\kappa(t)= 2^{-1}(t^{2}+2t)e^{-t}$, and $\varphi(t)=\Phi^{\prime}(t)$, the standard Gaussian density.

The aim of this short note is to extend above results to the case that $\rho_{ni}$ satisfies \eqref{eq1.4}. Notation such as $\widetilde{H}_{\lambda=0,\infty}(x,y)$ given by \eqref{minL} and $H(x,y)$ given by \eqref{eq1.5} will be used throughout this paper.

Contents of this paper are organized as follows. The main results are given in Section \ref{sec2} and Section \ref{sec3} presents the proofs.

\section{Main results}
In this section, we provide the main results. The first result is the joint limit distributions of minima and maxima as $\rho_{ni}$ is given by \eqref{eq1.4}.
\label{sec2}
\begin{theorem}
\label{th1}
Let the norming constant $b_{n}$ be given by \eqref{eq1.2}, under the assumption \eqref{eq1.4}, for every $(x_{i},y_{i})\in \R^{2},i=1,2$, with $\mathbf{v}_{n}$ and $\mathbf{u}_{n}$ given by \eqref{eq1.3}, we have
\begin{itemize}
\item[(i)] if \peng{$\lim_{n\to\infty}\max_{1\leq i \leq n} m(i/n)= 0$},
\begin{eqnarray*}
\lim_{n\to\infty}\P\left(\mathbf{M}_{n}\leq \mathbf{u}_{n},\mathbf{m}_{n}\leq\mathbf{v}_{n}\right)=H_{0}(x_{2},y_{2})\widetilde{H}_{0}(x_{1},y_{1});
\end{eqnarray*}

\item[(ii)] if \peng{$\lim_{n\to\infty} \min_{1\leq i \leq n} m(i/n)= \infty$},
\begin{eqnarray*}
\lim_{n\to\infty}\P\left(\mathbf{M}_{n}\leq \mathbf{u}_{n},\mathbf{m}_{n}\leq\mathbf{v}_{n}\right)=H_{\infty}(x_{2},y_{2})\widetilde{H}_{\infty}(x_{1},y_{1});
\end{eqnarray*}

\item[(iii)] if $ m(x)$ is a continuous positive function on [0,1],
\begin{eqnarray}
\label{add1}
\lim_{n\to\infty}\P\left(\mathbf{M}_{n}\leq \mathbf{u}_{n},\mathbf{m}_{n}\leq\mathbf{v}_{n}\right)=H\left(x_{2},y_{2}\right)\widetilde{H}(x_{1},y_{1});
\end{eqnarray}
where $H(x,y)$ is given by \eqref{eq1.5} and
\begin{equation}\label{min}
\widetilde{H}(x,y)=1-\Lambda(-x)-\Lambda(-y)+H(-x,-y).
\end{equation}
\end{itemize}
\end{theorem}

\begin{remark}
\label{remark2}
Note that $\rho_{ni}=1-\frac{m(i/n)}{\log n}$ implies $0\le m(i/n)\le 2\log n$. Examples of $m(x)$ satisfying the conditions mentioned in Theorem \ref{th1} are given as follows.
\begin{itemize}
\item[(1).] Case (i), let
\[m(i/n)=\left\{ {{\begin{array}{*{20}c}
 {i/n,\qquad \quad\qquad\qquad  i\in [1, n^{1/2}],} \hfill \\
 {1/n \bigwedge i/n= 1/n,\qquad \mbox{otherwise}.} \hfill \\
\end{array} }}\right. \]
Hence, $\lim_{n\to\infty}\max m(i/n)=0$. Actually, $m(i/n)$ given here satisfies the condition given by Theorem \ref{th3}.
\item[(2).] Case (ii), let
\[m(i/n)=\left\{ {{\begin{array}{*{20}c}
 {\log (n/i),\qquad \quad\qquad\qquad  i\in [1, n^{1/2}],} \hfill \\
 {(\log n) \bigvee \log (n/i)=\log n,\qquad  \mbox{otherwise}.} \hfill \\
\end{array} }}\right. \]
So,  $\lim_{n\to\infty}\min m(i/n)= \infty$. Note that $m(i/n)$ given here satisfies the condition given by Theorem \ref{th4}.
\item[(3).] For case (iii), let $m(x)=x+1$ with $0\le x\le 1$.
\end{itemize}
\end{remark}

\begin{remark}
\label{remark2}
Theorem \ref{th1} shows that $\mathbf{M}_{n}$ and $\mathbf{m}_{n}$ are \peng{asymptotically independent}, similar to the results of Hashorva and Weng (2013) and Liao and Peng (2014b) with  $\rho_{ni}$ satisfying \eqref{HRC}.
Similar results for  univariate weak dependent stationary case were proved by Davis (1979).
\end{remark}

\begin{corollary}
\label{corollary1}
Under the conditions of Theorem 1, we have
\[
\lim_{n\to\infty}\P\left(\max|\xi_{ni}|\leq b_{n}+\frac{x+\log 2}{b_{n}},\max|\eta_{ni}|\leq b_{n}+\frac{y+\log 2}{b_{n}}\right)
=H(x,y).
\]
\end{corollary}

Followings are convergence rates of joint distributions of maxima and minima to its ultimate extreme value distribution. There are three cases: $\lim_{n\to\infty}\max m(i/n)= 0$; $\lim_{n\to\infty}\min m(i/n)=\infty$; and
$m(t)$ is continuous and monotone on $[0,1]$.
\begin{theorem}
\label{th2}
Suppose that \eqref{eq1.4} holds. Further assume that $m(t)$ is continuous and monotone on $[0,1]$, then
\begin{eqnarray}
\label{eq2.1}
&&
\lim_{n\to\infty} \frac{2 \log n}{\log\log n}\left( \P\left(\mathbf{M}_{n}\leq\mathbf{u}_{n},\mathbf{m}_{n}\leq\mathbf{v}_{n}\right)-H(x_{2},y_{2})\widetilde{H}(x_{1},y_{1})\right)  \nonumber \\
&=&
H(x_{2},y_{2})\widetilde{H}(x_{1},y_{1}) e^{-x_{2}}\int_{0}^{1} m(t) \varphi\left(\sqrt{m(t)}+\frac{y_{2}-x_{2}}{2\sqrt{m(t)}}\right)dt \nonumber \\
&&+
H(x_{2},y_{2})H(-x_{1},-y_{1}) e^{x_{1}}\int_{0}^{1} m(t) \varphi\left(\sqrt{m(t)}+\frac{x_{1}-y_{1}}{2\sqrt{m(t)}}\right)dt,
\end{eqnarray}
where $\mathbf{v}_{n}$ and $\mathbf{u}_{n}$ are those given by \eqref{eq1.3}.
\end{theorem}

As $\lim_{n\to\infty}\max_{1\leq i\leq n} m(i/n)= 0$, we have the following result.
\begin{theorem}
\label{th3}
Assume that $\lim_{n\to\infty}(\log n)^{4} \max_{1\leq i\leq n} m(i/n)=0$ holds. With $\mathbf{v}_{n}$ and $\mathbf{u}_{n}$ given by \eqref{eq1.3} we have
\begin{eqnarray}
\label{eq2.2}
&&
\lim_{n\to\infty} (4\log n) \left(\P\left(\mathbf{M}_{n}\leq\mathbf{u}_{n},\mathbf{m}_{n}\leq\mathbf{v}_{n}\right)-H_{0}(x_{2},y_{2})\widetilde{H}_{0}(x_{1},y_{1})\right) \nonumber \\
&=&
\left\{\left[\min(x_{2},y_{2})\right]^{2}+2\min(x_{2},y_{2})\right\} e^{-\min(x_{2},y_{2})}H_{0}(x_{2},y_{2})\widetilde{H}_{0}(x_{1},y_{1}) \nonumber\\
&&
-\left\{\left[x_{1}^{2}-2x_{1}\right]e^{x_{1}}\Lambda(-x_{1})+\left[y_{1}^{2}-2y_{1}\right]e^{y_{1}}\Lambda(-y_{1})\right\}H_{0}(x_{2},y_{2}) \nonumber\\
&&+\left\{\left[\max(x_{1},y_{1})\right]^{2}-2\max(x_{1},y_{1})\right\} e^{\max(x_{1},y_{1})}H_{0}(x_{2},y_{2})H_{0}(-x_{1},-y_{1}) .
\end{eqnarray}
\end{theorem}

For the remainder case $\lim_{n\to\infty}\min_{1\leq i \leq n} m(i/n)=\infty$, we have the following result.
\begin{theorem}
\label{th4}
Let $\mathbf{v}_{n}$ and $\mathbf{u}_{n}$ be given by \eqref{eq1.3}. If $\lim_{n\to\infty}\left(\log\log n\right)/\min_{1\leq i \leq n} m(i/n)=0$, then
\begin{eqnarray}
\label{eq2.3}
&&
\lim_{n\to\infty} (4\log n) \left(\P\left(\mathbf{M}_{n}\leq\mathbf{u}_{n},\mathbf{m}_{n}\leq\mathbf{v}_{n}\right)-H_{\infty}(x_{2},y_{2})\widetilde{H}_{\infty}(x_{1},y_{1})\right)  \nonumber\\
&=&
\left\{\left[x_{2}^{2}+2x_{2}\right] e^{-x_{2}}+\left[y_{2}^{2}+2y_{2}\right] e^{-y_{2}}\right\}H_{\infty}(x_{2},y_{2})\widetilde{H}_{\infty}(x_{1},y_{1}) \nonumber\\
&&-\left\{\left[x_{1}^{2}-2x_{1}\right] e^{x_{1}}\Lambda(-x_{1}) + \left[y_{1}^{2}-2y_{1}\right] e^{y_{1}}\Lambda(-y_{1})\right\} H_{\infty}(x_{2},y_{2}) \nonumber\\
&&+\left\{\left[x_{1}^{2}-2x_{1}\right] e^{x_{1}}+\left[y_{1}^{2}-2y_{1}\right] e^{y_{1}}\right\}H_{\infty}(x_{2},y_{2})H_{\infty}(-x_{1},-y_{1}).
\end{eqnarray}
\end{theorem}

\section{Proofs}
\label{sec3}

\noindent {\bf Proof of Theorem \ref{th1}.}  Let $F_{i}(x,y)$ denote the bivariate Gaussian distribution with correlation coefficient $\rho_{ni}$, i.e., the joint distribution of $(\xi_{ni}, \eta_{ni}), 1\leq i \leq n$, and let $u_{n}(x)=b_{n}+x/b_{n}$, $v_{n}(x)=-b_{n}+x/b_{n}$ for notational simplicity.

We first show (iii). It follows from Theorem 1 in Liao and Peng (2016) that
\begin{eqnarray}
\label{eq3.1}
&&
-\sum^{n}_{i=1}\left(1-F_{i}\left(u_{n}(x),u_{n}(y)\right)\right)
\nonumber\\
&&\to
-e^{-y}\int_{0}^{1}\Phi\left(\sqrt{m(t)}+\frac{x-y}{2\sqrt{m(t)}}\right)dt
-e^{-x}\int_{0}^{1}\Phi\left(\sqrt{m(t)}+\frac{y-x}{2\sqrt{m(t)}}\right)dt ,n\to\infty.
\end{eqnarray}
By using \eqref{eq1.2} and assumption \eqref{eq1.4}, we have
\begin{eqnarray}
\label{eq3.3}
\lim_{n\to\infty}\frac{u_{n}(x_{2})-\rho_{ni}v_{n}(z)}{\sqrt{1-\rho_{ni}^{2}}}=\infty.
\end{eqnarray}
It follows from \eqref{eq3.3} and the dominated convergence theorem that
\begin{eqnarray}
\label{eq3.4}
&&
\sum^{n}_{i=1} \P\left(\xi_{ni}>u_{n}(x_{2}),\eta_{ni}\leq v_{n}(y_{1})\right)  \nonumber\\
&=&
\sum^{n}_{i=1}\int^{v_{n}(y_{1})}_{-\infty}\left(1-\Phi\left(\frac{u_{n}(x_{2})-\rho_{ni}z}{\sqrt{1-\rho_{ni}^{2}}} \right)\right)\varphi(z) dz \nonumber \\
&=&
\Big(1+b_{n}^{-2}+O(b_{n}^{-4})\Big)\int_{-\infty}^{y_{1}} n^{-1} \sum^{n}_{i=1}
\left[1-\Phi\left(\frac{u_{n}(x_{2})-\rho_{ni}v_{n}(z)}{\sqrt{1-\rho_{ni}^{2}}}\right)\right] \exp{(z-\frac{z^{2}}{2b_{n}^{2}})} dz  \nonumber\\
&\to&
0
\end{eqnarray}
as $n\to\infty$ due to the fact that $\rho_{ni}\to 1$ uniformly for all $i$.

Similarly we have
\begin{eqnarray}
\label{eq3.5}
\lim_{n\to\infty}\sum^{n}_{i=1} \P\left(\xi_{ni}\leq v_{n}(x_{1}),\eta_{ni}>u_{n}(y_{2})\right)=0.
\end{eqnarray}

It follows from \eqref{eq3.1}, \eqref{eq3.4}, \eqref{eq3.5} and $(-\xi_{ni},-\eta_{ni})\stackrel{d}{=}(\xi_{ni},\eta_{ni})$ that
\begin{eqnarray}
&&
\sum^{n}_{i=1}\log \P\left(v_{n}(x_{1})<\xi_{ni}\leq u_{n}(x_{2}),v_{n}(y_{1})<\eta_{ni}\leq u_{n}(y_{2})\right)  \nonumber\\
&=&
-\sum^{n}_{i=1}\Big[1- \P\left(v_{n}(x_{1})<\xi_{ni}\leq u_{n}(x_{2}),v_{n}(y_{1})<\eta_{ni}\leq u_{n}(y_{2})\right) \nonumber\\ &&+\frac{1}{2}
\left(1- \P\left(v_{n}(x_{1})<\xi_{ni}\leq u_{n}(x_{2}),v_{n}(y_{1})<\eta_{ni}\leq u_{n}(y_{2})\right)\right)^{2}(1+o(1))\Big]  \nonumber\\
&=&
-\sum^{n}_{i=1}\Big[1-F_{i}\left(u_{n}(x_{2}),u_{n}(y_{2})\right)+1-F_{i}\left(u_{n}(-x_{1}),u_{n}(-y_{1})\right)\nonumber\\
&&-\P\left(\xi_{ni}>u_{n}(x_{2}),\eta_{ni}\leq v_{n}(y_{1})\right) -
\P\left(\xi_{ni}\leq v_{n}(x_{1}),\eta_{ni}>u_{n}(y_{2})\right) \nonumber\\ &&+ \frac{1}{2}\left(1- \P\left(v_{n}(x_{1})<\xi_{ni}\leq u_{n}(x_{2}),v_{n}(y_{1})<\eta_{ni}\leq u_{n}(y_{2})\right)\right)^{2}(1+o(1))\Big]  \nonumber\\
&\to&
-e^{-y_{2}}\int_{0}^{1}\Phi\left(\sqrt{m(t)}+\frac{x_{2}-y_{2}}{2\sqrt{m(t)}}\right)dt
-e^{-x_{2}}\int_{0}^{1}\Phi\left(\sqrt{m(t)}+\frac{y_{2}-x_{2}}{2\sqrt{m(t)}}\right)dt  \nonumber\\
&&-e^{y_{1}}\int_{0}^{1}\Phi\left(\sqrt{m(t)}+\frac{y_{1}-x_{1}}{2\sqrt{m(t)}}\right)dt
-e^{x_{1}}\int_{0}^{1}\Phi\left(\sqrt{m(t)}+\frac{x_{1}-y_{1}}{2\sqrt{m(t)}}\right)dt \label{eq3.6}
\end{eqnarray}
as $n\to\infty$, hence \eqref{eq3.6} implies
\begin{eqnarray}
\label{eq3.7}
\lim_{n\to\infty} \P(\mathbf{v}_{n}< \mathbf{m}_{n}\leq \mathbf{M}_{n} \leq \mathbf{u}_{n} )=H(x_{2},y_{2})H(-x_{1},-y_{1}).
\end{eqnarray}

Noting that
\begin{eqnarray}
\label{eq3.8}
&&
\sum^{n}_{i=1}\log \P\left(v_{n}(x_{1})<\xi_{ni}\leq u_{n}(x_{2}),\eta_{ni}\leq u_{n}(y_{2})\right)  \nonumber\\
&=&
-\sum^{n}_{i=1}\left(1-F_{i}\left(u_{n}(x_{2}),u_{n}(y_{2})\right)\right)-n\Phi\left(v_{n}(x_{1})\right)+
\sum^{n}_{i=1} \P\left(\xi_{ni}\leq v_{n}(x_{1}),\eta_{ni}> u_{n}(y_{2})\right)  \nonumber\\
&&-
\frac{1}{2}\sum^{n}_{i=1}\left(1-\P\left(v_{n}(x_{1})<\xi_{ni}\leq u_{n}(x_{2}),\eta_{ni}\leq u_{n}(y_{2})\right)\right)^{2}(1+o(1))
\end{eqnarray}
and
\begin{eqnarray}
\label{eq3.12}
\lim_{n\to\infty}n\left(1-\Phi(u_{n}(x))\right)=e^{-x},\quad\quad \lim_{n\to\infty}n\Phi(v_{n}(x))=e^{x},
\end{eqnarray}
so by \eqref{eq3.1}, \eqref{eq3.5}, \eqref{eq3.8} and \eqref{eq3.12}, we have
\begin{eqnarray}
\label{eq3.9}
\lim_{n\to\infty} \P( \mathbf{M}_{n} \leq \mathbf{u}_{n}, \mathbf{m}_{n1}>\mathbf {v}_{n}(x_{1}) )=H(x_{2},y_{2})\Lambda(-x_{1}).
\end{eqnarray}

Similarly,
\begin{eqnarray}
\label{eq3.10}
\lim_{n\to\infty}\P( \mathbf{M}_{n} \leq \mathbf{u}_{n}, \mathbf{m}_{n2}>\mathbf {v}_{n}(y_{1}) )=H(x_{2},y_{2})\Lambda(-y_{1}).
\end{eqnarray}

Finally, combing \eqref{eq1.5}, \eqref{eq3.7}, \eqref{eq3.9} and \eqref{eq3.10}, we can get
\begin{eqnarray*}
\lim_{n\to\infty}\P( \mathbf{M}_{n} \leq \mathbf{u}_{n}, \mathbf{m}_{n}\leq
\mathbf{v}_{n}) &=& \lim_{n\to\infty} \Big(\P(\mathbf{M}_{n} \leq \mathbf{u}_{n})
- \P(\mathbf{M}_{n} \leq \mathbf{u}_{n}, m_{n1}>v_{n}(x_{1}))\nonumber\\
& -&  \P(\mathbf{M}_{n} \leq \mathbf{u}_{n}, m_{n2}>v_{n}(y_{1})) +
\P(\mathbf{v}_{n} <\mathbf{m}_{n}\leq \mathbf{M}_{n} \leq
\mathbf{u}_{n})\Big)  \\ \nonumber
&=&
H(x_{2},y_{2})\widetilde{H}(x_{1},y_{1}).
\end{eqnarray*}
The proof of case (iii) is complete.

Now we return to case (i). By arguments similar to the proof of case (iii), we have
\begin{eqnarray}
\label{eq3.11}
\lim_{n\to\infty}\sum^{n}_{i=1}\left(1-F_{i}(u_{n}(x),u_{n}(y))\right) = e^{-\min(x,y)}.
\end{eqnarray}

By using (C.2) in Piterbarg (1996), we have
\[
\sum^{n}_{i=1} \P\left(\xi_{ni}>u_{n}(x_{2}),\eta_{ni}\leq v_{n}(y_{1})\right) \le
n\Phi(v_{n}(y_{1}))\Phi(-u_{n}(x_{2})) \to 0
\]
as $n\to\infty$ due to \eqref{eq3.12} and $\rho_{ni}>0$.

Combining \eqref{eq1.5}, \eqref{add1}, \eqref{eq3.8} and \eqref{eq3.11}, we finish the proof of case (i).
The proof of case (ii) is similar. Details are omitted here.

The proof of Theorem \ref{th1} is complete.
\qed

\noindent {\bf Proof of Theorem \ref{th2}.} By Theorem 2 of Liao and Peng (2016), we have
\begin{eqnarray}
\label{eq4.1}
&&
-\sum^{n}_{i=1}\left(1-F_{i}\left(u_{n}(x),u_{n}(y)\right)\right)+e^{-x}+\int^{1}_{0}\int^{\infty}_{y}
\Phi\left(\sqrt{m(t)}+\frac{x-z}{2\sqrt{m(t)}}\right)e^{-z}dz dt   \nonumber\\
&&\sim
\frac{\log\log n}{2\log n} e^{-x} \int^{1}_{0}\sqrt{m(t)}\varphi\left(\sqrt{m(t)}+\frac{y-x}{2\sqrt{m(t)}}\right) dt,\;\; n\to\infty
\end{eqnarray}
and Nair (1981) showed that
\begin{eqnarray}
\label{eq4.2}
\lim_{n\to\infty} b_{n}^{2}\left(n\Phi(v_{n}(x))-e^{x}\right)=-\frac{x^{2}-2x}{2}e^{x}.
\end{eqnarray}

Noting that by using Mills' inequality, for sufficient large $n$ and fixed $z$ we have
\begin{eqnarray}
\label{eq4.4}
&&
b_{n}^{2}\left[1-\Phi\left(\frac{u_{n}(x_{2})-\rho_{ni}v_{n}(z)}{\sqrt{1-\rho_{ni}^{2}}}\right)\right] \nonumber\\
&<&
b_{n}^{2}\frac{\sqrt{1-\rho_{ni}^{2}}}{\sqrt{2\pi}\left[(1+\rho_{ni})b_{n}+\frac{1}{b_{n}}(x_{2}-\rho_{ni}z)\right]}
\exp\left(-\frac{(1+\rho_{ni})^{2}b_{n}^{2}+2(1+\rho_{ni})(x_{2}-\rho_{ni}z)}{2(1-\rho_{ni}^{2})}\right)  \nonumber\\
&<&
\frac{b_{n}^{2}\exp\left(-\frac{(1+\rho_{ni})b_{n}^{2}+2(x_{2}-\rho_{ni}z)}{2(1-\rho_{ni})}-\frac{1}{2}\log\frac{b_{n}^{2}}{1-\rho_{ni}}\right)}
{\sqrt{2\pi}\left(1+\frac{1}{b_{n}^{2}(1+\rho_{ni})}(x_{2}-\rho_{ni}z)\right)}  \nonumber\\
&<&
\frac{\exp\left(-\frac{(1+\rho_{ni})b_{n}^{2}+2(x_{2}-\rho_{ni}z)}{2(1-\rho_{ni})}+\log b_{n}\sqrt{1-\rho_{ni}}\right)}
{\sqrt{2\pi}\left(1+\frac{1}{2 b_{n}^{2}}(x_{2}-\rho_{ni}z)\right)}  \nonumber\\
&\to&
0
\end{eqnarray}
as $n\to\infty$ since assumptions of Theorem \ref{th2} implies that $\rho_{ni}\to 1$ and $ b_{n}\sqrt{1-\rho_{ni}}$ is bounded uniformly for all $i$.
Hence, by using \eqref{eq4.4} we have
\begin{eqnarray}
\label{eq4.5}
&&
 b_{n}^{2} \sum^{n}_{i=1} \P\left(\xi_{ni}>u_{n}(x_{2}),\eta_{ni}\leq v_{n}(y_{1})\right)  \nonumber\\
&=&
\left(1+b_{n}^{-2}+O(b_{n}^{-4})\right)\int_{-\infty}^{y_{1}} n^{-1} \sum^{n}_{i=1} b_{n}^{2}
\left[1-\Phi\left(\frac{u_{n}(x_{2})-\rho_{ni}v_{n}(z)}{\sqrt{1-\rho_{ni}^{2}}}\right)\right] \exp{(z-\frac{z^{2}}{2b_{n}^{2}})} dz  \nonumber\\
&\to&
0
\end{eqnarray}
as $n\to\infty$. Similarly,
\begin{eqnarray}
\label{eq4.6}
\lim_{n\to\infty} b_{n}^{2} \sum^{n}_{i=1} \P\left(\xi_{ni}\leq v_{n}(x_{1}),\eta_{ni}>u_{n}(y_{2})\right) = 0.
\end{eqnarray}
Hence, by \eqref{eq4.1}, \eqref{eq4.2}, \eqref{eq4.5} and \eqref{eq4.6}, we  have
\begin{eqnarray}
\label{eq4.7}
&&
\P\left(\mathbf{M}_{n}\leq \mathbf{u}_{n},{m}_{n1}>v_{n}(x_{1})\right)-H(x_{2},y_{2})\Lambda(-x_{1})  \nonumber\\
&=&
H(x_{2},y_{2})\Lambda(-x_{1})\Big[\exp\Big(\sum^{n}_{i=1}\log \P(v_{n}(x_{1})<\xi_{ni}\leq u_{n}(x_{2}),\eta_{ni}\leq u_{n}(y_{2}))\nonumber\\ &&+e^{-x_{2}}+
\int^{1}_{0}\int^{\infty}_{y_{2}}\Phi\left(\sqrt{m(t)}+\frac{x_{2}-z}{2\sqrt{m(t)}}\right)e^{-z} dz dt +e^{x_{1}} \Big)-1\Big]  \nonumber\\
&=&
H(x_{2},y_{2})\Lambda(-x_{1})\Big[\exp \Big(-\sum^{n}_{i=1}(1-F_{i}(u_n(x_{2}),u_n(y_{2}))) -n\Phi(v_{n}(x_{1}))\nonumber \\&&+\sum^{n}_{i=1} \P\left(\xi_{ni}\leq v_{n}(x_{1}),\eta_{ni}>u_{n}(y_{2}))\right)-\sum^{n}_{i=1}\frac{1}{2}\left(1-\P(v_{n}(x_{1})<\xi_{ni}\leq u_{n}(x_{2}),\eta_{ni}\leq u_{n}(y_{2}))\right)^{2}(1+o(1)) \nonumber\\ &&+e^{-x_{2}}+
\int^{1}_{0}\int^{\infty}_{y_{2}}\Phi\left(\sqrt{m(t)}+\frac{x_{2}-z}{2\sqrt{m(t)}}\right)e^{-z} dz dt+e^{x_{1}} \Big)-1\Big]  \nonumber \\
&=&
H(x_{2},y_{2})\Lambda(-x_{1})\Big[-\sum^{n}_{i=1}(1-F_{i}(u_n(x_{2}),u_n(y_{2})))+e^{-x_{2}}+
\int^{1}_{0}\int^{\infty}_{y_{2}}\Phi\left(\sqrt{m(t)}+\frac{x_{2}-z}{2\sqrt{m(t)}}\right)e^{-z} dz dt  \nonumber\\
&&+ e^{x_{1}}-n\Phi(v_{n}(x_{1}))+\sum^{n}_{i=1} \P\left(\xi_{ni}\leq v_{n}(x_{1}),\eta_{ni}>u_{n}(y_{2}))\right)  \nonumber \\ &&-\sum^{n}_{i=1}\frac{1}{2}\left(1-\P(v_{n}(x_{1})<\xi_{ni}\leq u_{n}(x_{2}),\eta_{ni}\leq u_{n}(y_{2}))\right)^{2}(1+o(1))\Big](1+o(1)) \nonumber\\
&\sim& H(x_{2},y_{2})\Lambda(-x_{1})\Big(\frac{\log\log n}{2\log n} e^{-x_{2}}\int^{1}_{0}\sqrt{m(t)}\varphi\left(\sqrt{m(t)}+\frac{y_{2}-x_{2}}{2\sqrt{m(t)}}\right)dt
+\frac{1}{2\log n}\frac{x_{1}^{2}-2 x_{1}}{2}e^{x_{1}}\Big).
\end{eqnarray}

Similarly,
\begin{eqnarray}
\label{eq4.8}
&&
\P\left(\mathbf{M}_{n}\leq \mathbf{u}_{n},{m}_{n2}>v_{n}(y_{1})\right)-H(x_{2},y_{2})\Lambda(-y_{1})  \nonumber\\
&&\sim
H(x_{2},y_{2})\Lambda(-y_{1})\Big(\frac{\log\log n}{2\log n} e^{-x_{2}}\int^{1}_{0}\sqrt{m(t)}\varphi\left(\sqrt{m(t)}+\frac{y_{2}-x_{2}}{2\sqrt{m(t)}}\right)dt \nonumber\\
&&\quad+\frac{1}{2\log n}\frac{y_{1}^{2}-2 y_{1}}{2}e^{y_{1}}\Big).
\end{eqnarray}

It follows from \eqref{eq3.6} and \eqref{eq4.1} that
\begin{eqnarray}
\label{eq4.9}
&&
\P\left(\mathbf{v}_{n}<\mathbf{m}_{n}\leq \mathbf{M}_{n}\leq\mathbf{u}_{n}\right)-H(x_{2},y_{2})H(-x_{1},-y_{1})  \nonumber\\
&=&
H(x_{2},y_{2})H(-x_{1},-y_{1})\Big[\exp\Big(\log\prod_{i=1}^{n}\P(v_{n}(x_{1})<\xi_{ni}\leq u_{n}(x_{2}),v_{n}(y_{1})<\eta_{ni}\leq u_{n}(y_{2}))) \nonumber\\ &&-\log H(x_{2},y_{2})H(-x_{1},-y_{1})\Big )-1\Big]  \nonumber\\
&=&
H(x_{2},y_{2})H(-x_{1},-y_{1})\Big[\exp\Big(\sum^{n}_{i=1}\log\P(v_{n}(x_{1})<\xi_{ni}\leq u_{n}(x_{2}),v_{n}(y_{1})<\eta_{ni}\leq u_{n}(y_{2})) \nonumber\\&&+
e^{-x_{2}} +\int^{1}_{0}\int^{\infty}_{y_{2}}\Phi\left(\sqrt{m(t)}+\frac{x_{2}-z}{2\sqrt{m(t)}}\right)e^{-z} dz dt \nonumber\\
&&+ e^{x_{1}}+\int^{1}_{0}\int^{\infty}_{-y_{1}}\Phi\left(\sqrt{m(t)}+\frac{-x_{1}-z}{2\sqrt{m(t)}}\right)e^{-z} dz dt\Big)-1\Big]  \nonumber\\
&=&
H(x_{2},y_{2})H(-x_{1},-y_{1})\Big[-\sum^{n}_{i=1}(1-F_{i}(u_{n}(x_{2}),u_{n}(y_{2})))-\sum^{n}_{i=1}(1-F_{i}(u_{n}(-x_{1}),u_{n}(-y_{1})))  \nonumber\\
&&+\sum^{n}_{i=1}\P(u_{n}(x_{2})<\xi_{ni},\eta_{ni}\leq v_{n}(y_{1}))+\sum^{n}_{i=1}\P(\xi_{ni}\leq v_{n}(x_{1}),u_{n}(y_{2})<\eta_{ni})+e^{-x_{2}}  \nonumber\\
&&+\int^{1}_{0}\int^{\infty}_{y_{2}}\Phi\left(\sqrt{m(t)}+\frac{x_{2}-z}{2\sqrt{m(t)}}\right)e^{-z} dzdt+e^{x_{1}}+\int^{1}_{0}\int^{\infty}_{-y_{1}}\Phi\left(\sqrt{m(t)}+\frac{-x_{1}-z}{2\sqrt{m(t)}}\right)e^{-z} dz dt  \nonumber\\
&&-\sum^{n}_{i=1}\frac{1}{2}\left(1-\P(v_{n}(x_{1})<\xi_{ni}\leq u_{n}(x_{2}),v_{n}(y_{1})<\eta_{ni}\leq u_{n}(y_{2}))\right)^{2}(1+o(1))\Big](1+o(1))  \nonumber\\
&\sim&
\frac{\log\log n}{2\log n}H(x_{2},y_{2})H(-x_{1},-y_{1})\Big[e^{-x_{2}}\int^{1}_{0}\sqrt{m(t)}\varphi\left(\sqrt{m(t)}+\frac{y_{2}-x_{2}}{2\sqrt{m(t)}}\right)dt \nonumber\\
&&+
e^{x_{1}}\int^{1}_{0}\sqrt{m(t)}\varphi\left(\sqrt{m(t)}+\frac{x_{1}-y_{1}}{2\sqrt{m(t)}}\right)dt\Big].
\end{eqnarray}
Combing Theorem 2 of Liao and Peng (2016) with \eqref{eq4.7}-\eqref{eq4.9}, we can get \eqref{eq2.1}. The proof is complete.
\qed

\noindent {\bf Proof of Theorem \ref{th3}.}
By arguments similar  to that of  \eqref{eq4.5}-\eqref{eq4.9}, we may derive the following facts:
\begin{eqnarray*}
&&
\P\left(\mathbf{M}_{n}\leq \mathbf{u}_{n},{m}_{n1}>v_{n}(x_{1})\right)-H_{0}(x_{2},y_{2})\Lambda(-x_{1})  \nonumber\\
&&\sim\frac{1}{4\log n} H_{0}(x_{2},y_{2})\Lambda(-x_{1})\Big\{\left[\left(\min(x_{2},y_{2})\right)^{2}+2\min(x_{2},y_{2})\right]e^{-\min(x_{2},y_{2})}+\left[x_{1}^{2}-2x_{1}\right]e^{x_{1}}\Big\};
\end{eqnarray*}
\begin{eqnarray*}
&&
\P\left(\mathbf{M}_{n}\leq \mathbf{u}_{n},{m}_{n2}>v_{n}(y_{1})\right)-H_{0}(x_{2},y_{2})\Lambda(-y_{1})  \\
&&\sim \frac{1}{4\log n}H_{0}(x_{2},y_{2})\Lambda(-y_{1})\Big\{\left[\left(\min(x_{2},y_{2})\right)^{2}+2\min(x_{2},y_{2})\right]e^{-\min(x_{2},y_{2})}+\left[y_{1}^{2}-2y_{1}\right]e^{y_{1}}\Big\}
\end{eqnarray*}
and
\begin{eqnarray*}
&&\P\left(\mathbf{v}_{n}<\mathbf{m}_{n}\leq \mathbf{M}_{n}\leq\mathbf{u}_{n}\right)-H_{0}(x_{2},y_{2})H_{0}(-x_{1},-y_{1}) \\
&\sim&
\frac{1}{4\log n}H_{0}(x_{2},y_{2})H_{0}(-x_{1},-y_{1})\Big\{\left[\left(\min(x_{2},y_{2})\right)^{2}+2\min(x_{2},y_{2})\right]e^{-\min(x_{2},y_{2})} \\
&&\qquad\qquad+\left[\left(\max(x_{1},y_{1})\right)^{2}-2\max(x_{1},y_{1})\right]e^{\max(x_{1},y_{1})}\Big\}.
\end{eqnarray*}
Combining those facts with Theorem 3 in Liao and Peng (2016), we can get \eqref{eq2.2}. The proof is complete.
\qed

\noindent {\bf Proof of Theorem \ref{th4}.}
By arguments similar to \eqref{eq4.5}-\eqref{eq4.9},
we may derive the following facts:
\begin{eqnarray*}
&&
\P\left(\mathbf{M}_{n}\leq \mathbf{u}_{n},{m}_{n1}>v_{n}(x_{1})\right)-H_{\infty}(x_{2},y_{2})\Lambda(-x_{1})  \nonumber\\
&&\sim \frac{1}{4\log n}H_{\infty}(x_{2},y_{2})\Lambda(-x_{1})\Big\{\left[x_{2}^{2}+2x_{2}\right]e^{-x_{2}}+\left[y_{2}^{2}+2y_{2}\right]e^{-y_{2}}+\left[x_{1}^{2}-2x_{1}\right]e^{x_{1}}\Big\};
\end{eqnarray*}
\begin{eqnarray*}
&&
\P\left(\mathbf{M}_{n}\leq \mathbf{u}_{n},{m}_{n2}>v_{n}(y_{1})\right)-H_{\infty}(x_{2},y_{2})\Lambda(-y_{1})  \\
&&\sim \frac{1}{4\log n}H_{\infty}(x_{2},y_{2})\Lambda(-y_{1})\Big\{\left[x_{2}^{2}+2x_{2}\right]e^{-x_{2}}+\left[y_{2}^{2}+2y_{2}\right]e^{-y_{2}}+\left[y_{1}^{2}-2y_{1}\right]e^{y_{1}}\Big\}
\end{eqnarray*}
and
\begin{eqnarray*}
&&\P\left(\mathbf{v}_{n}<\mathbf{m}_{n}\leq \mathbf{M}_{n}\leq\mathbf{u}_{n}\right)-H_{\infty}(x_{2},y_{2})H_{\infty}(-x_{1},-y_{1}) \\
&\sim&
\frac{1}{4\log n}H_{\infty}(x_{2},y_{2})H_{\infty}(-x_{1},-y_{1})\Big\{\left[x_{2}^{2}+2x_{2}\right]e^{-x_{2}}+\left[y_{2}^{2}+2y_{2}\right]e^{-y_{2}}
+\left[x_{1}^{2}-2x_{1}\right]e^{x_{1}}+\left[y_{1}^{2}-2y_{1}\right]e^{y_{1}}\Big\}.
\end{eqnarray*}
Combining those facts with Theorem 4 in Liao and Peng (2016),
we can get \eqref{eq2.3}, and the proof is complete.
\qed


\begin{thebibliography}{999}

\bibitem{}
Engelke, S., Kabluchko, Z. and Schlather, M. (2015). Maxima of independent, non-identically distributed Gaussian vectors. {\it Bernoulli}, {\bf 1}, 38--61.

\bibitem{}
Davis, R. (1979). Maxima and minima of statinary sequences. {\it The Annals of Probability}, {\bf 3}, 453--460.


\bibitem{}
de Haan, L. and Peng, L. (1997). Rates of convergence for bivariate extremes. {\it Journal of Multivariate Analysis}, {\bf 61}, 195--230.


\bibitem{}
de Haan, L., and Resnick, S. I. (1996). Second-order regular variation and rates of convergence in extreme-value theory.
{\it The Annals of Probability}, {\bf 24}, 97--124.

\bibitem{}
D\c{e}bicki, K., Hashorva, E. and Ji, L. (2014). Gaussian approximation of perturbed chi-square risk. {\it Statistics and Its Interface}, {\bf 7}, 363--373.

\bibitem{}
French, J. P. and Davis, R. (2013). The asymptotic distribution of the maxima of a Gaussian random field on a lattice. {\it Extremes}, {\bf 16}, 1--26.

\bibitem{}
Frick, M. and Riess, R.-D. (2013) Expansions and penultimate distributions of maxima of bivariate normal random vectors. {\it Statisitcs and Probability Letters}, {\bf 83}, 2563--2568.

\bibitem{}
Hall, P. (1979). On the rate of convergence of normal extremes. {\it Journal of Applied Probability}, {\bf 16}, 433--439.

\bibitem{}
Hashorva, E. (2005). Elliptical triangular arrays in the max-domain
of attraction of H\"{u}sler-Reiss distribution. {\it Statistics and
Probability Letters}, {\bf 72}, 125--135.

\bibitem{}
Hashorva, E. (2006). On the max-domain of attractions of bivariate
elliptical arrays. {\it Extremes}, {\bf 8}, 225--233.

\bibitem{}
Hashorva, E. (2013). Minima and maxima of elliptical arrays and
spherical processes. {\it Bernoulli}, {\bf 3}, 886--904.

\bibitem{}
Hashorva, E., Kabluchko, Z. and W\"{u}bker, A. (2012). Extremes of independent chi-square random vectors. {\it Extremes}, {\bf 15}, 35--42.

\bibitem{}
Hashorva, E., Nadarajah, S. and Pogani, T.K. (2014). Extremes of perturbed bivariate rayleigh risks. {\it REVSTAT--Statistical Journal}, {\bf 2}, 157--168

\bibitem{}
Hashorva, E., Peng, Z. and Weng, Z. (2016). Higher-order expansions
of distributions of maxima in a H\"{u}sler-Reiss model. {\it
Methodology and Computing in Applied Probability}, {\bf 18}, 181-196.

\bibitem{}
Hashorva, E. and Weng, Z. (2013). Limit laws for extremes of dependent stationary Gaussian arrays. {\it Statistics and Probability Letters}, {\bf 83}, 320--330.

\bibitem{}
Hsing, T., H\"{u}sler, J. and Reiss, R.-D. (1996). The extremes of a triangular array of normal random variables. {\it The Annals of Applied Probability}, {\bf 2}, 671--686.

\bibitem{}
H\"usler, J. and Reiss, R.-D. (1989). Maxima of normal random
vectors: between independence and complete dependence. {\it Statistics and Probability Letters}, {\bf 7}, 283--286.

\bibitem{}
Kabluchko, Z., de Haan, L. and Schlatter, M. (2009). Stationary max-stable fields associated to negative
definite functions. {\it The Annals of Probability}, {\bf 37}, 2042--2065.


\bibitem{}
Liao, X. and Peng, Z. (2014a). Convergence rate of maxima of
bivariate Gaussian arrays to the H\"{u}sler-Reiss distribution. {\it
Statistics and Its Interface}, {\bf 7}, 351--362.

\bibitem{}
Liao, X. and Peng, Z. (2014b). Asymptotics for the maxima and minima
of H\"{u}sler-Reiss bivariate Gaussian arrays. {\it Extremes}, {\bf 18},
1--14.


\bibitem{}
Liao, X. and Peng, Z. (2016). Asymptotics and statistical
inferences on independent and non-identically distributed bivariate
Gaussian triangular arrays. http://arxiv.org/abs/1505.03431v3

\bibitem{}
Liao, X., Peng, L., Peng, Z. and Zheng, Y. (2016). Dynamic bivariate normal copula.
{\it Science China Mathematics}, {\bf 5}, 955--976.

\bibitem{}
Nair, K. A. (1981). Asymptotic distribution and moments of normal
extremes. {\it The Annals of Probability}, {\bf 9}, 150--153.

\bibitem{}
Piterbarg, V. I. (1996). {\it Asymptotic Methods in the Theory of Gaussian Processes and Fields}. In: Translations of Mathematical Monographs, Volume 148. American Mathematical Society, Providence, Rhode Island.

\end{thebibliography}
\end{document}